\RequirePackage[leqno]{amsmath}
\documentclass[leqno]{amsart}
\usepackage{amsmath}
\makeatletter
\usepackage[utf8]{inputenc}     
\usepackage{multicol}
\usepackage{enumitem}
\usepackage{color, graphics}
\usepackage[demo]{graphicx}
\usepackage{efbox,graphicx}
\efboxsetup{linecolor=gray,linewidth=0.8pt}
\usepackage{floatrow}
\usepackage{pgfplots}
\usepackage{graphicx, multicol, latexsym, amsmath, amssymb}
\usepackage{pgfplots}
\usepackage{tikz}
\usepackage{amsfonts}
\usepackage{enumitem}
\setlist[itemize]{topsep=0pt,after=\vspace{1.5\baselineskip}}
 \usepackage[colorlinks=true]{hyperref}
\usepackage{empheq}
\usepackage[T1]{fontenc}
\usepackage{tabularx}
\usepackage[leftcaption]{sidecap}
\sidecaptionvpos{figure}{m}
\usepackage{tikz}
\usepackage{subfigure}
\usepackage[margin=1in]{geometry}
\usetikzlibrary{backgrounds,automata}

\newcount\rc@count
\let\rc@clearconstantlist\empty
\newcommand\rc@clearconstant[1]{\global\expandafter\let\csname rc@const@#1\endcsname\undefined}
\newcommand\resetconstants[1]{%
    \def\rc@constname{#1}
    \global\rc@count=1\relax 
    \bgroup 
        \let\\\rc@clearconstant 
        \rc@clearconstantlist
        \global\let\rc@clearconstantlist\empty 
    \egroup
}
\newcommand\const[1]{%
    \@ifundefined{rc@const@#1}{%
        \expandafter\xdef\csname rc@const@#1\endcsname{%
           \noexpand\rc@useconst{\rc@constname}{\the\rc@count}%
        }%
        \g@addto@macro\rc@clearconstantlist{\\{\mathrm{#1}}}%
        \global\advance\rc@count1\relax
    }{}%
    \csname rc@const@#1\endcsname
}
\newcommand\rc@useconst[2]{{#1}\textsubscript{#2}}
\setlist[itemize]{noitemsep, topsep=0pt}

\def\R{\mathbb R} \def\N{\mathbb N}

\def\R{\mathbb R} \def\N{\mathbb N} 
\def\TM{T_\text{max}} 
\def
\@cite
#1#2{[#1\if@tempswa, #2\fi]}
\makeatother
\newtheorem{theorem}{Theorem}[section]

\newtheorem{lemma}[theorem]{Lemma}

\newtheorem{remark}{Remark}

\newcounter{cnstcnt}

\title[Nonlocal effects on a chemotaxis consumption model]{
Global dynamics of chemotaxis-consumption systems with oppositely acting nonlocal terms
}
\author[Rafael D\'iaz Fuentes, Fatma Gamze D\"{u}zg\"{u}n, Silvia Frassu and Giuseppe Viglialoro]{}
\keywords{Chemotaxis, Global existence, Nonlocal growth terms, Boundedness. \\
\textit{$^*$Corresponding author}: silvia.frassu@unica.it}

\begin{document}
\maketitle
\centerline{\scshape{\scshape{Rafael D\'iaz Fuentes$^{\flat}$, Fatma Gamze D\"{u}zg\"{u}n$^{\flat}$, Silvia Frassu$^{\flat,*}$ \and Giuseppe Viglialoro$^{\flat}$}}}
\medskip
{
\medskip
\centerline{$^\flat$Dipartimento di Matematica e Informatica}
\centerline{Universit\`{a} degli Studi di Cagliari}
\centerline{Via Ospedale 72, 09124. Cagliari (Italy)}
\medskip
}
\bigskip
\begin{multicols}{2}
\tableofcontents
\end{multicols}
\begin{abstract}
In this paper, we investigate a chemotaxis system featuring a class of external source terms that incorporate both local and nonlocal growth as well as dampening effects. The model describes the evolution of a cell density \( u \), which migrates in response to a chemical signal \( v \), within an impermeable habitat. Mathematically, this leads to the study of the following initial-boundary value problem
\begin{equation}\label{problem_abstract}
\tag{$\Diamond$}
\begin{cases}
u_t = \Delta u - \chi \nabla \cdot (u \nabla v) +f(u)& \text{in } \Omega \times (0, T_{\max}), \\
 v_t = \Delta v - v u & \text{in } \Omega \times (0, T_{\max}), \\
u_\nu  = v_\nu  = 0 & \text{on } \partial\Omega \times (0, T_{\max}), \\
u(x, 0) = u_0(x) \geq 0,\quad v(x,0) = v_0(x) \geq 0 & \text{for } x \in \bar{\Omega},
\end{cases}
\end{equation}
where $x\in \Omega$, and $\Omega$ is a bounded and smooth domain of $\R^n$ ($n \in \mathbb{N}$) with boundary $\partial\Omega$, oriented by the outward unit normal vector $\nu$, and $\chi>0$. Moreover for $a, b > 0$, $\alpha, \beta, \gamma \geq 1$ the function $f(u)$ is assumed to take one of the following forms:
\[
\text{either } \; a u^\alpha - b u^\beta\int_\Omega u^\gamma \textrm{ or } -a u^\alpha + b u^\beta\int_\Omega u^\gamma.\] 
Finally, $\TM\in (0,\infty]$ indicates the maximal instant of time up to which the corresponding solutions exist.
 
We provide conditions ensuring the absence of aggregation phenomena over time. Specifically, we show that the maximal existence time satisfies $T_{\max} = \infty$, and both $u$ and $v$ remain uniformly bounded for all time, in the following situations: 
\begin{enumerate}[label=($\Diamond_{\arabic*}$)]
\item  \label{1Abstract} for $f(u)=a u^\alpha - b u^\beta\int_\Omega u^\gamma$, under the assumption $\beta\geq \alpha$, whenever either $1 \leq \alpha < 2$ and $\beta+\gamma > \frac{n}{2}+2$, or $\alpha \geq 2$ and $\beta +\gamma> \frac{n}{2}(\alpha-1)+\alpha$;
\item \label{2Abstract} for $f(u)=-a u^\alpha + b u^\beta\int_\Omega u^\gamma$, in these scenarios: $\alpha>2$ and $\beta+\gamma<\alpha$, $\alpha>2$, $\beta+\gamma=\alpha$ and $a>b|\Omega|$, and finally $\alpha=2$, $\beta=\gamma=1$ and $a>b|\Omega|+ \mathcal{C}_P \| \chi v_0\|_{L^\infty(\Omega)}$ for some $\mathcal{C}_P=\mathcal{C}_P(n,\Omega)>0$.
\end{enumerate}
This work builds upon previous studies that established global existence and boundedness results for variants of system~\eqref{problem_abstract}, in which the second equation takes the form 
$\tau v_t = \Delta v - v + u$ 
(cf.~\cite{Bian2018} for the case $\tau = 0$ and~\cite{ChiyoEtAl2024Nonlocal} for $\tau = 1$). 
In those settings, the source term $f$ typically has a more \textit{ad hoc} structure, namely 
$f(u) = a u^{\alpha} - b u^{\alpha} \int_{\Omega} u^{\beta}$. 
Based on this framework, the present study focuses on distinguishing the different dynamical behaviors exhibited by the systems in term of the reaction described in \ref{1Abstract} and \ref{2Abstract}. 
In the first case, the model structure ensures the control of the total mass $\int_{\Omega} u$ over time, whereas in the second such control is no longer guaranteed, leading to substantially different qualitative dynamics. 
\end{abstract}
\resetconstants{c}
\section{Introduction and selected known results}\label{IntroSection}
\subsection{An overview on Keller--Segel models}
The mathematical modeling of chemotaxis was first introduced by Keller and Segel in the 1970s~\cite{Keller1970,Keller-Segel1970}. They proposed a system of partial differential equations governing the spatiotemporal evolution of a biological population \( u = u(x,t) \) (which may represent cells, microorganisms, or bacteria) in response to a chemical substance \( v = v(x,t) \), commonly termed a chemoattractant. The system is considered within a bounded and sufficiently smooth domain \( \Omega \subset \mathbb{R}^n \) over a time interval \( (0, T_{\max}) \); naturally, $x \in \Omega$ and $t \in (0, T_{\max})$, where $T_{\max}$ represents the maximal existence time of the solution, that is, the supremum of all times for which $u$ and $v$ remain well defined.

More precisely, the most representative equation for the classical Keller--Segel model reads as follows
\begin{equation}\label{problemOriginalKS}
u_t = \Delta u - \chi \nabla \cdot (u \nabla v) \quad \text{in } \Omega \times (0,T_{\max}),
\end{equation}
and it describes the population dynamics incorporating linear diffusion and directed movement driven by chemotaxis, with \( \chi > 0 \) denoting the chemotactic sensitivity. The term \( -\chi \nabla \cdot (u \nabla v) \) represents the flux of the population moving up the gradient of the chemoattractant concentration.

This equation is coupled with one of the following models describing the chemoattractant dynamics: a production mechanism
\begin{equation}\label{SecondEqKS}
v_t = \Delta v - v + u \quad \text{in } \Omega \times (0,T_{\max}),
\end{equation}
or an absorption mechanism
\begin{equation}\label{SecondEqKS-cons}
v_t = \Delta v - v u \quad \text{in } \Omega \times (0,T_{\max}),
\end{equation}
where the production/absorption terms are linear in \( u \) and idealize respectively the release or consumption of the chemical signal by the population. To ensure that the system of equations is well-posed, it is necessary to supplement it with appropriate boundary conditions and initial data for $u$ and $v$. In particular, homogeneous Neumann boundary conditions indicate that the process takes place in an impermeable domain, while $
u(x,0) = u_0(x)$ and $v(x,0) = v_0(x)$
specify the initial distributions of the cell density and the chemical signal, respectively.

In models incorporating chemoattractant production, the interplay between diffusive dispersal and chemotactic aggregation gives rise to complex spatiotemporal patterns. The outcome of this competition depends sensitively on the spatial dimension \( n \), the initial total mass $
m := \int_{\Omega} u_0(x) \, dx,$
and the chemotactic sensitivity parameter \( \chi \). Under certain conditions, this interaction may lead to blow-up phenomena at finite time ($\TM$ finite), characterized by the formation of singularities and unbounded population densities.

In the one-dimensional setting, global-in-time ($\TM=\infty$) existence and uniform boundedness of solutions are well established. However, as the spatial dimension increases, the dynamics become substantially more intricate. Specifically, for \( n \geq 3 \), solutions may blow up in finite time even when the initial mass is relatively small (see the survey \cite[$\S$3]{BellomoEtAlM3AS} and references therein). 

The two-dimensional case has been extensively studied due to the emergence of the \emph{critical mass phenomenon}, a threshold effect whereby the global behavior of solutions--whether they exist globally in time or exhibit finite-time blow-up--depends critically on the size of the initial total mass (see, e.g., \cite{Jager1992,Nagai}).

\subsection{An overview of Keller--Segel models with production and logistic effects}
Incorporating logistic-type source terms into the classical Keller--Segel system governed by equations \eqref{problemOriginalKS} and \eqref{SecondEqKS} is known to enhance the prospects for global existence by introducing a damping effect.
\subsubsection{Classical logistics}
A typical example tied to what now said, is the modified equation
\begin{equation*}\label{KS}
u_t = \Delta u - \chi \nabla \cdot (u \nabla v) + au - b u^\beta \quad \text{in } \Omega \times (0,\TM), \textrm{ where } a,b>0 \textrm{ and } \beta\geq 1,
\end{equation*}
which has been widely studied as a regularized variant of the original chemotaxis model. While the presence of the damping term indeed improves the regularity properties of solutions, global boundedness has been rigorously established only under specific conditions. In particular, for the case \( \beta = 2 \), global existence holds provided that the coefficient \( b \) is sufficiently large (see \cite{TelloWinkParEl, W0}). However, blow-up may still occur when \( \beta \) is close to 1; see \cite{WinDespiteLogistic, FuestCriticalNoDEA, Winkler_ZAMP-FiniteTimeLowDimension}.
\subsubsection{Nonlocal logistics}
In parallel, attention has shifted to source terms involving nonlocal damping mechanisms of the form
\begin{equation}\label{NonLocalSources}
f(u) = a u^\alpha - b u^\alpha \int_\Omega u^\beta, \quad \text{with } \quad \alpha, \beta \ge 1.
\end{equation}
Such expression models population growth regulated by a death rate that depends on the total biomass, and it can be viewed as nonlocal generalizations of logistic damping. The presence of the integral term introduces a spatially global feedback that may significantly alter the long-time dynamics of solutions. 

A central analytical question is whether such nonlocal damping is strong enough to prevent finite-time blow-up and guarantee boundedness of solutions; that is, for the equation 
\begin{equation*}\label{KS-NonLocal}
u_t = \Delta u - \chi \nabla \cdot (u \nabla v) + a u^\alpha - b u^\alpha \int_\Omega u^\beta,
\end{equation*}
can global boundedness arise solely from the nonlocal term, or are further constraints required? Studies incorporating nonlinear diffusion, chemotactic sensitivity, signal production, or multiple chemical signals establish that the nonlocal term can prevent blow-up under suitable conditions \cite{NegreanuTelloNonLocalPara-Ell, Bian2018, TaoFang-NonlinearNonlocal, Tello2021}. (Also in $\mathbb{R}^n$, nonlocal damping with Poisson-type chemoattractant interactions has a regularizing effect: \cite{BianEtAlWholeSpaceNonLocal,ChenWangDocMat,CarrilloWangAcAppMat,Li-Viglialoro-DIE-Nonlocal}.)
\subsubsection{Gradient-dependent logistics}
Furthermore, finite-time blow-up can be prevented by incorporating a logistic-type source $h=h(u, |\nabla u|)$ into the population cells' equation $u_t = \Delta u - \chi \nabla \cdot (u \nabla v)$. Typical forms include
$
h(u, |\nabla u|) = a u^\alpha - b u^\beta - c |\nabla u|^\gamma$,
or alternatively, $
h(u, |\nabla u|) = a u^\alpha - b u^\alpha \int_\Omega u^\beta - c |\nabla u|^\gamma$,
with positive constants \( a,b,c \) and exponents \( \alpha, \beta, \gamma \geq 1\). While gradient-dependent nonlinearities arising from external sources often interpreted biologically as mechanisms of \emph{accidental death} have been extensively studied within single-equation models (see, e.g., \cite{QS} and references therein), their influence within chemotactic systems has been lately garnering significant attention. Notably, recent investigations such as \cite{IshidaLankeitVigliloro-Gradient, LiEtAl2024gradientnonlinearitiesprevent, BaghaeiEtAlJDE_2026, DiazEtAl_2025-AAM} have begun to elucidate their role in biologically and ecologically relevant contexts.
\begin{remark}[Consumption models with classical and gradient-dependent logistics]\label{RemarkProdutcion}
Above, we have summarized the results we consider most relevant to the present study, focusing on Keller–Segel type models in which the chemical signal(s) are produced by the same cell population. In parallel, for consumption-type systems, namely those in which the signal equation is of the form \eqref{SecondEqKS-cons}, there also exist investigations incorporating classical logistic terms (cf. \cite{LankeitWang_DCDS-B,BaghaeiKhelgati_MMAS_2017,BaghaeiCRMA_2023}), as well as studies addressing gradient-dependent logistic mechanisms (\cite{Columbu_JMAA_2025}). However, to the best of our knowledge, no study has yet addressed consumption models involving nonlocal logistic terms. As we shall emphasize later, this gap in the literature represents one of the main motivations of the present work.
\end{remark}
\section{Presentation of the model and motivations; claim of the main result}
In accordance with Remark \ref{RemarkProdutcion}, and as it will be discussed in detail below, our analysis is devoted to consumption-type Keller–Segel systems endowed with nonlocal source terms. Furthermore, we shall consider logistic components in which the nonlocality constitutes, on the one hand, a natural generalization of that appearing in \eqref{NonLocalSources}, and, on the other hand, may influence both proliferative and dissipative dynamics within the population.
\subsection{Basic description of the model}
In this paper, we study the initial-boundary value problem
\begin{equation}\label{problem}
\begin{cases}
u_t=\Delta u-\chi\nabla \cdot (u\nabla v)+f(u) &\text{in } \Omega \times (0, \TM),\\
v_t=\Delta v-vu &\text{in } \Omega \times (0, \TM),\\
u_\nu=v_\nu=0 &\text{on } \partial\Omega \times (0, \TM),\\
u(x, 0)=u_0(x),\quad v(x,0)=v_0(x) &x \in \bar{\Omega},
\end{cases}
\end{equation}
where $x\in \Omega \subset \mathbb{R}^n$ ($n \in \mathbb{N}$), $\Omega$ is a bounded domain with smooth boundary $\partial \Omega$, and the parameter $\chi$ is positive. The initial data $u_0(x), v_0(x)$ are assumed to be nonnegative and sufficiently regular and for $a, b > 0$, $\alpha, \beta, \gamma \geq 1$ the function $f(u)$ may take one of the following forms:
\begin{equation}\label{NonLocalSourcesHEreinUsed}
\text{either } \; a u^\alpha - b u^\beta\int_\Omega u^\gamma \textrm{ or } -a u^\alpha + b u^\beta\int_\Omega u^\gamma.
\end{equation}
As anticipated, and reiterated here for the sake of clarity, the notation \((\cdot)_\nu\) denotes the outward normal derivative on \(\partial \Omega\), while \(\TM\) represents the maximal existence time of the solution.

In line with the previous sections, the first equation models the population density $u$ via diffusion, chemotaxis, and local growth or decay, while the second governs the chemoattractant $v$ through diffusion and consumption. Homogeneous Neumann conditions enforce no flux on $\partial \Omega$.
\subsubsection{Interpretation of \( f(u) \) and its role in the investigation}
The expressions of $f(u)$ in \eqref{NonLocalSourcesHEreinUsed} describe nonlocal damping or growth mechanisms, where the evolution of $u(x,t)$ depends not only on the local value but also on the global quantity $\int_\Omega u^\gamma$, representing the total population mass raised to a power and introducing spatially averaged feedback.

Precisely, in the first scenario the term \( a u^\alpha \) models local population growth, typically driven by reproduction or proliferation processes. The accompanying nonlocal inhibitory term \( -b u^\beta \int_\Omega u^\gamma \) serves as a saturation mechanism, reflecting effects such as overcrowding or competition for limited resources that depend on the overall population size. Exactly in the opposite direction, the second formulation features a leading local decay term \( -a u^\alpha \), which may be interpreted as natural mortality or dissipation, while the nonlocal term \( b u^\beta \int_\Omega u^\gamma \) can induce growth once the total population surpasses a certain threshold. 

Despite differences in the source term, both forms fit a unified framework for analyzing boundedness and global existence. However, from a mathematical point of view, the choice of logistic term in \eqref{NonLocalSourcesHEreinUsed} is crucial and determines the technical approach. In the first case, uniform-in-time boundedness of the mass, $\int_\Omega u$, holds for all $t \in (0, T_{\max})$, while in the second case this property is lost. Accordingly, the analysis must be adapted depending on the form of $f(u)$ (see Remark \ref{RemarkMassa} and $\S$\ref{SectionAPrioriEstimates} for details).
\subsection{Formulation of the main Theorem}
In order to present our result, here we assume 
\begin{equation}\label{reglocal}
\begin{cases}
    \textrm{for some } \sigma \in(0,1) \textrm{ and } n\in \N, \Omega \subset \R^n \textrm{ is a bounded domain of class } C^{2+\sigma}, \textrm{ with boundary } \partial \Omega, \\
    u_0, v_0: \bar{\Omega}  \rightarrow \R^+_0 , \textrm{ with } u_0,  v_0 \in C^{2+\sigma}(\bar\Omega)\; \textrm{such that} \; u_{0\nu}=0\; \textrm{and} \; v_{0\nu}=0\;  \textrm{on }\partial \Omega. 
   \end{cases}
\end{equation}
Subsequently, we will prove the following  
\begin{theorem}\label{theoremlocal}
Let the hypotheses in \eqref{reglocal} be fulfilled and $\chi, b>0$ and $\alpha,\gamma \geq 1$.  Then problem \eqref{problem} admits a unique solution
\begin{equation*}
(u,v)\in C^{2+\sigma,1+\frac{\sigma}{2}}( \Bar{\Omega} \times [0, \infty))\times C^{2+\sigma,1+\frac{\sigma}{2}}(\Bar{\Omega} \times [0, \infty))
\end{equation*}
with $0\leq u,v \in L^\infty(\Omega \times (0,\infty))$, in these circumstances: 
\begin{enumerate}[label=\bf\arabic*., ref=\bf\arabic*]
\item For $f(u)=a u^\alpha - b u^\beta\int_\Omega u^\gamma$, with any $a>0$ and $\beta\geq \alpha$  whenever 
\begin{enumerate} [label=(\alph*), ref=\theenumi(\alph*)]
\item  \label{a} $\alpha \geq 2$ and $\beta + \gamma> \frac{n}{2}(\alpha-1)+\alpha$;
\item   \label{b} $1 \leq \alpha < 2$ and $\beta + \gamma > \frac{n}{2}+2$;
\end{enumerate}
\item For $f(u)=-a u^\alpha + b u^\beta\int_\Omega u^\gamma$, whenever $\beta \geq 1$ and
\begin{enumerate}[label=(\alph*), ref=\theenumi(\alph*)]
\item \label{(a)} $\alpha>2$, $\beta+\gamma<\alpha$, with any $a>0$;  
\item \label{(b)} $\alpha>2$, $\beta+\gamma=\alpha$, with $a>b|\Omega|$; 
\item \label{(c)} $\alpha=2$, $\beta=\gamma=1$, with $a>b|\Omega|+ \mathcal{C}_P \| \chi v_0\|_{L^\infty(\Omega)}$ for some $\mathcal{C}_P=\mathcal{C}_P(n,\Omega)>0$.
\end{enumerate} 
\end{enumerate}
\end{theorem} 
\begin{remark}
In the absence of a logistic term, i.e., when $f(u) = 0$, our model \eqref{problem} reduce to pure consumption mechanisms, which remain an interesting subject of study, particularly regarding the existence of blow-up solutions. In this context, in \cite{TaoConsumption2011} and \cite{BaghaeiKhelgati_CR-NoLogistic} the authors obtained boundedness under suitable smallness assumptions on $\chi \|v_0\|_{L^\infty(\Omega)}$, while in \cite{LankeitWang_DCDS-B} it is established that, in the presence of a classical logistic term of the form $b u - a u^2$, boundedness is guaranteed provided $a > \varphi(\chi \|v_0\|_{L^\infty(\Omega)},n)$, for some continuous expression of $\varphi(\xi,n)$, 
$\xi >0.$ Notably, the limit scenario obtained for $b=0$ in case \ref{(c)} of Theorem \ref{theoremlocal} qualitatively recovers \cite[Theorem 1.1]{LankeitWang_DCDS-B}.
\end{remark}
\section{Preliminary results}
The upcoming auxiliary results will be used in the proof of our main estimates. 
Moreover, we shall adopt the following conventions:
\begin{itemize}
\item[$\lozenge$] all constants \( c_i \), \( i = 1, 2, \ldots \), appearing below are tacitly assumed to be positive;
\item[$\lozenge$] the symbol \( \epsilon \) will denote an arbitrary positive real number; 
multiplication by another constant or addition of another arbitrary constant will not be performed, 
and for convenience, the final result will also be denoted by \( \epsilon \).
\end{itemize}
\begin{lemma}\label{GagliardoIneqLemma}
Let $\Omega$ complying with the hypotheses in \eqref{reglocal}, and for $n\geq 3$ let  
\begin{equation*}\label{def_of_p}
p:=\frac{2n}{n-2}. 
\end{equation*}
Additionally, let $q, r$ satisfy $1 \le r<q<p$ and $\frac{q}{r}<\frac{2}{r}+1-\frac{2}{p}$. Then for all $\epsilon>0$ there exists $C_0=C_0(\epsilon)>0$ such that for all $\varphi \in H^1(\Omega) \cap L^r(\Omega)$, 
\begin{equation}\label{Bian}
\|\varphi\|_{L^q(\Omega)}^q \le C_0 \|\varphi\|_{L^r(\Omega)}^\delta
+\epsilon \|\nabla \varphi\|_{L^2(\Omega)}^2+\|\varphi\|_{L^2(\Omega)}^2, 
\end{equation}
where 
\begin{equation*}
\lambda:=\frac{\frac{1}{r}-\frac{1}{q}}{\frac{1}{r}-\frac{1}{p}} \in (0,1), \quad \delta:=\frac{2(1-\lambda)q}{2-\lambda q}.
\end{equation*}
The same conclusion holds for $n\in\{1,2\}$ whenever
$q, r$ fulfill, respectively,  $1 \le r<q$ and $\frac{q}{r}<\frac{2}{r}+2$ and $1 \le r<q$ and $\frac{q}{r}<\frac{2}{r}+1$.
\begin{proof}
The proof can be found in \cite[Lemma 3.1]{ChiyoEtAl2024Nonlocal}.
\end{proof}
\end{lemma}
\begin{lemma}\label{LemmaODI-Comparison}
Let\/ $T>0$ and $\phi:(0,T)\times \R^+_0\rightarrow \R$. If $0\leq y\in C^0([0,T))\cap  C^1((0,T))$ is such that 
\begin{equation*}
y'\leq \phi(t,y)\quad \textrm{for all } t \in (0,T), 
\end{equation*}
and there is $y_1>0$ with the property that whenever $y>y_1$ for some $t\in (0,T)$ one has that $\phi(t,y)\leq 0$, then
\begin{equation*}
y\leq \max\{y_1,y(0)\}\quad \textrm{on } (0,T).
\end{equation*}
\begin{proof}
For the proof we refer the reader to \cite[Lemma 3.3]{ChiyoEtAl2024Nonlocal}.
\end{proof}
\end{lemma}
\section{Local existence, parabolic regularity results and globality}\label{SectionLocalExis}
Let us establish existence and uniqueness of local classical solutions to problem \eqref{problem}. We refrain from presenting the full proof, as it is routine in this framework, and instead we emphasize the main steps.
\begin{lemma}\label{theoremExistence}
Let the hypotheses in \eqref{reglocal} be fulfilled.  Additionally, let $\chi, a, b>0$ and $\alpha,\beta,\gamma \geq 1$ and $f(u)$ as in \eqref{NonLocalSourcesHEreinUsed}.  Then 
there exist $\TM\in (0,\infty]$ and a unique couple of functions $(u,v$), with 
\begin{equation*}
(u,v)\in C^{2+\sigma,1+\frac{\sigma}{2}}( \Bar{\Omega} \times [0, \TM))\times C^{2+\sigma,1+\frac{\sigma}{2}}( \Bar{\Omega} \times [0, \TM)),
\end{equation*}
solving problem \eqref{problem}, and complying with 
\begin{equation}\label{boundednessv}
u\geq 0 \quad \textrm{and}\quad 0\leq v \leq ||v_0||_{L^\infty(\Omega)} \quad \textrm{in}\quad \bar{\Omega}\times [0,\TM). 
\end{equation}
Finally,
 \begin{equation}\label{dictomyCriteC2+del} 
 \text{if} \quad \TM<\infty \quad \text{then} \quad \limsup_{t \to \TM} \left(\|u(\cdot,t)\|_{L^\infty(\Omega)}\right)=\infty.
 \end{equation}
\begin{proof}
Classical fixed-point arguments in H\"{o}lder spaces ensure the local well-posedness of the problem, establishing the existence and uniqueness of a classical solution $(u, v)$ on \([0,T]\) for some sufficiently small \(T>0\). By virtue of the maximum principle, the structural condition \(f(0)=0\), together with the homogeneous Neumann boundary condition and the assumption \(u_0 \geq 0\) on \(\bar{\Omega}\), guarantees the nonnegativity of \(u\) for all \(t \in [0,T]\); consequently, the same principle yields the second relation in \eqref{boundednessv}. Standard prolongation techniques ensure the existence of a maximal time of existence \(\TM \in (0,\infty]\); furthermore, if \(\TM < \infty\), this time corresponds to a finite-time blow-up, in the sense of \eqref{dictomyCriteC2+del}. For more details we indicate \cite[Proposition~4]{Bian2018}.
\end{proof}
\end{lemma}
From this point onward \((u,v)\) will denote the unique local solution to problem~\eqref{problem}, defined on \(\Omega \times (0,\TM)\) and guaranteed by Lemma~\ref{theoremExistence}.
\begin{lemma}\label{lem:MaxReg} 
For every $q>1$ there exists $C_P=C_P(n,\Omega,q)>0$ for which $v$ is such that
\begin{equation*}\label{tau1}
\int_0^t e^s \int_\Omega |\Delta v(\cdot,s)|^q\,ds \leq C_{P}2^{q-1}\|v_0\|^{q}_{L^\infty(\Omega)}\int_0^t e^s \int_\Omega u(\cdot,s)^q\,ds+\const{gs}e^t \quad \text{for all } t\in(0,\TM).
\end{equation*}
Additionally, if $u\in L^\infty((0,\TM);L^k(\Omega))$ for some $k>n$, then
\begin{equation}\label{tau1extension}
v \in L^\infty((0,\TM);W^{1,\infty}(\Omega)).
\end{equation}
\begin{proof}  
Since $v$ solves $v_t = \Delta v - v u$ in $\Omega \times (0,\TM)$, it equivalently solves $v_t = \Delta v - v + v(1-u)$ in $\Omega \times (0,\TM)$. By invoking maximal Sobolev regularity results for this latter equation (see \cite{hieber_pruess} or \cite[Theorem~2.3]{giga_sohr}) and following the methodology outlined in \cite{IshidaYokotaDCDS-SViaMSR,IshidaYokotaJDE-2012SmallData,SenbaSuzukiAAN2006}, one obtains
\begin{equation}\label{uuoo}
\int_0^t e^s \int_\Omega |\Delta v(\cdot,s)|^q \, ds \le C_P  
\Bigl(\|v_0\|_{1 - 1/q, q}  + \int_0^t e^s \int_\Omega |v(1-u(\cdot,s))|^q \, ds \Bigr), \quad \textit{for all } t \in (0,\TM),
\end{equation}
where $C_P= C_P(n,\Omega,q)$ and $\|\cdot\|_{1-1/q,q}$ denotes the real interpolation norm.
Exploiting the upper bound in \eqref{boundednessv}, and $(A+B)^q\leq 2^{q-1}(A^q+B^q)$, with $A,B\geq 0$ and $q\geq 1$ (and such an inequality might be employed in other places of the paper without explicit mention), we have that $|v(1-u(\cdot,s))|^q \le 2^{q-1}\|v_0\|^q_{L^\infty(\Omega)} (1+u(\cdot,s)^q)$, 
and by substituting this into \eqref{uuoo} yields the desired conclusion.

As for the final claim, we first note that for $k > n$ it holds that 
$- \tfrac{1}{2} - \tfrac{n}{2k} > -1$.  Making use of the standard $L^k$–$L^\infty$ estimate for the Neumann heat semigroup (see, for instance, \cite[Lemma~1.3]{win_aggregation_vs}), we infer that
\begin{align*}
\|\nabla v(\cdot, t)\|_{L^\infty(\Omega)}
&\le \|\nabla e^{t\Delta} v(\cdot,0)\|_{L^\infty(\Omega)}
   + \int_0^t \|\nabla e^{(t-s)\Delta} \big(u(\cdot,s) v(\cdot,s)\big)\|_{L^\infty(\Omega)} \, ds \\
&\le \|\nabla v_0\|_{L^\infty(\Omega)}
   + \const{c_1} \int_0^t 
     \left( 1 + (t-s)^{-\frac{1}{2} - \frac{n}{2k}} \right)
     e^{-\mu_1 (t-s)}
     \|u(\cdot,s) v(\cdot,s)\|_{L^k(\Omega)} \, ds,
\end{align*}
for all $t \in (0,\TM)$, where $\mu_1$ denotes the first nonzero eigenvalue of 
$-\Delta$ in $\Omega$ under homogeneous Neumann boundary conditions.  
Since, again, $v$ is bounded, the assumption on $u$ ensures that 
$\|u(\cdot,s) v(\cdot,s)\|_{L^k(\Omega)}$ remains uniformly controlled for 
$s \in (0,\TM)$, we may conclude that $
\|\nabla v(\cdot, t)\|_{L^\infty(\Omega)} \le \const{c_2}$ 
for all $t \in (0,\TM)$,  and thus the claim follows.
\end{proof}
\end{lemma}
\subsection{A boundedness criterion}
As specified, when relation~\eqref{dictomyCriteC2+del} holds the solution is said to exhibit blow-up at some finite time \(\TM\) in the \(L^{\infty}(\Omega)\)-norm.  The next result rules out this blow-up scenario and ensures globality of solutions.  
\begin{lemma}[A boundedness criterion] \label{ExtensionLemma} There exists $k_0=k_0(n,\alpha,\beta,\gamma)\geq n+2$ such that $u$ whenever  $u\in L^\infty((0,\TM);L^{k}(\Omega))$ for $k>k_0$ then we have $\TM=\infty$ and in particular $u\in L^\infty((0,\infty);L^\infty(\Omega))$.
\begin{proof}
Once it has been established that for all $k \geq k_0$ one has the inclusion $
u \nabla v \in L^{\infty}\!\big((0, \TM); L^{q_1}(\Omega)\big)$ 
(see ~\eqref{tau1extension}), with $q_1 > n + 2$, 
the uniform-in-time boundedness of $u$ in $L^{\infty}(\Omega)$ for all $t \in (0,\TM)$ 
follows by appropriately adjusting $k_0$ 
also in term of $\alpha,\beta,\gamma$ so as to satisfy all the structural assumptions of~\cite[Lemma~A.1]{TaoWinkParaPara}. Subsequently, by invoking the dichotomy criterion in \eqref{dictomyCriteC2+del} of Lemma \ref{theoremExistence}, we necessarily have $\TM=\infty$.
\end{proof}
\end{lemma} 
\begin{remark}[On the role of the value $k_0$ in Lemma \ref{ExtensionLemma}] 
Within the framework of our model, the expressions of $f(u)$ given in~\eqref{NonLocalSourcesHEreinUsed} effectively play the role of the source term $g(x,t)$ appearing in~\cite[Lemma~A.1]{TaoWinkParaPara}. 
Building upon this correspondence, we emphasize the following considerations, which will be relevant for the subsequent analysis:
\begin{enumerate}
\item The parameter $k_0$ should be regarded merely as a lower bound which, for certain combinations of $\alpha$, $\beta$, and $\gamma$, may suffice to guarantee the boundedness of the solution, although in other parameter regimes it may need to be enlarged. 
In this sense, $k_0$ is to be understood as a quantity which, while in principle computable, is neither optimal nor uniquely determined, since its precise value depends on the interplay among $n$, $\alpha$, $\beta$, and $\gamma$. Accordingly, when $f(u)$ is taken in the forms of~\eqref{NonLocalSourcesHEreinUsed}, no specific restrictions on $\alpha$, $\beta$, and $\gamma$ will be required, and the value of $k_0$ may be tacitly assumed to increase whenever necessary, so to ensures the mathematical consistency of the subsequent arguments. 
\item In contrast, when $f(u)$ is considered in the second form of~\eqref{NonLocalSourcesHEreinUsed} and  the corresponding analysis is confined to the case $\alpha = 2$ and $\beta = \gamma=1$, we can fix $k_0 = n + 2$.
\end{enumerate}
\end{remark}
\section{A priori estimates} \label{SectionAPrioriEstimates}
In this section, we derive a priori estimates for the local solution $(u,v)$ to model \eqref{problem}. 
These estimates yield bounds on the solution that depend essentially only on the prescribed initial data, and they represent a fundamental tool in establishing the uniform-in-time boundedness of the solution itself.
\begin{remark}\label{RemarkMassa}
We remark that, in the analysis that follows, a certain property of the solution is guaranteed in one case but not in the other, namely the uniform-in-time boundedness of the mass $\int_\Omega u$. Indeed, when $f(u)$ takes the first form in \eqref{NonLocalSourcesHEreinUsed}, the mass can be controlled; in the second case, this is no longer possible. As a result, the mathematical approach differs: in the first case, we will make use of the Gagliardo--Nirenberg inequality and the mentioned boundedness of the mass to handle, by means of the nonpositive nonlocal term of the source, contributions from the cross diffusion term and the nonnegative local part of the source itself. 
Conversely, in the second case, the nonnegative nonlocal term and the cross-diffusion contribution are controlled by the nonpositive local term of the logistic component, without it being previously required to ensure mass boundedness.
\end{remark}
In the next steps we will dedicate to derive uniform-in-time bounds of $u$ in some $L^k(\Omega)$-spaces, for $k>1.$
\subsection{The case $f(u)=a u^\alpha-bu^\beta \int_\Omega u^\gamma$}
\begin{lemma}[Boundedness of the mass]\label{localSol}
Let $\beta \geq \alpha \geq 1$ and $\gamma \geq 1$. Then there exists $m_0>0$ such that 
\begin{equation}\label{MassBounded}
\int_\Omega u(x,t)\,dx \leq m_0 \quad \textrm{for all } t \in (0,\TM).
\end{equation}
\begin{proof}
For $\beta>\alpha$, we integrate over $\Omega$ the first equation of problem \eqref{problem} so that by H\"{o}lder's inequalities,
\begin{equation*}\label{M1}
\begin{split}
y'(t):=\frac{d}{dt} \int_\Omega u &\leq a |\Omega|^{\frac{\beta-\alpha}{\beta}} \left(\int_\Omega u^\beta\right)^{\frac{\alpha}{\beta}} - b |\Omega|^{(1-\gamma)}\int_\Omega u^\beta \left(\int_\Omega u \right)^{\gamma}\\
&= \left(\int_\Omega u^\beta\right)^{\frac{\alpha}{\beta}} \left[a |\Omega|^{\frac{\beta-\alpha}{\beta}} - b |\Omega|^{(1-\gamma)} \left(\int_\Omega u^\beta\right)^{\frac{\beta-\alpha}{\beta}} \left(\int_\Omega u \right)^{\gamma}\right]\\ 
&\leq \left(\int_\Omega u^\beta\right)^{\frac{\alpha}{\beta}} \left[a |\Omega|^{\frac{\beta-\alpha}{\beta}} - b |\Omega|^{(1-\gamma)-\frac{(\beta-1)(\beta-\alpha)}{\beta}} 
\left(\int_\Omega u \right)^{\beta-\alpha+\gamma}\right] \quad \textrm{for all } t \in (0, \TM).
\end{split}
\end{equation*}
Now we apply Lemma \ref{LemmaODI-Comparison} with $T=\TM$, $\rho(t):= \left(\int_\Omega u^{\beta}\right)^{\frac{\alpha}{\beta}}\geq 0$ on $(0,\TM)$, 
$$\phi(t,y)=\rho(t)\left[a |\Omega|^{\frac{\beta-\alpha}{\beta}} - b |\Omega|^{(1-\gamma)-\frac{(\beta-1)(\beta-\alpha)}{\beta}} (y(t))^{\beta-\alpha+\gamma}\right],$$ $y_0=y(0)=\int_\Omega u_0$ and $y_1:=\left(\frac{a |\Omega|^{\frac{\beta-\alpha}{\beta}}}{b |\Omega|^{(1-\gamma)-\frac{(\beta-1)(\beta-\alpha)}{\beta}}}\right)^{\frac{1}{\beta-\alpha+\gamma}}$, so concluding with $m_0=\max\{y_0, y_1\}$.

For $\beta=\alpha$ and $\gamma \geq 1$, we directly obtain
\begin{equation*}\label{M1}
\begin{split}
y'(t):=\frac{d}{dt} \int_\Omega u &\leq \int_\Omega u^\alpha \left[a - b |\Omega|^{1-\gamma} \left(\int_\Omega u \right)^{\gamma}\right] \quad \textrm{for all } t \in (0, \TM),
\end{split}
\end{equation*}
so we achieve the claim by relying again in Lemma \ref{LemmaODI-Comparison}. 
\end{proof}
\end{lemma}
Let us start by considering the case $\alpha \ge 2$.
\begin{lemma}\label{CaseAlphaMaggiore2Lemma}
Assume that $\alpha, \beta, \gamma \ge 1$ satisfy that
\begin{equation}\label{alpha_large}
\beta \geq \alpha \geq 2 \quad \textrm{and} \quad \beta + \gamma>\frac{n(\alpha-1)}{2}+ \alpha.
\end{equation}
Then $u\in L^\infty((0,\TM);L^{k}(\Omega))$, for all $k>k_0$.
\begin{proof}
From the first equation in \eqref{problem}, an integration by parts yields for all $k>k_0$\begin{equation}\label{uk_estimate1}
\begin{split}
\frac{d}{dt}\int_\Omega u^k &=-k(k-1)\int_\Omega u^{k-2}|\nabla u|^2+k(k-1)\chi\int_\Omega u^{k-1}\nabla u\cdot \nabla v +ka\int_\Omega u^{k+\alpha-1}\\
&\quad -kb\left(\int_\Omega u^{k+\beta-1}\right)\left(\int_\Omega u^\gamma\right)\\
&=-\frac{4(k-1)}{k}\int_\Omega |\nabla u^{\frac{k}{2}}|^2-(k-1)\chi\int_\Omega u^k\Delta v +ka\int_\Omega u^{k+\alpha-1}\\
&\quad -kb\left(\int_\Omega u^{k+\beta-1}\right)\left(\int_\Omega u^\gamma\right) 
\quad  \mbox{on}\ (0, \TM). 
\end{split}
\end{equation}
Since $\alpha>1$, the Young inequality gives
\begin{equation}\label{uLap_v_estimate}
-(k-1) \chi \int_\Omega u^k\Delta v \le 
\int_\Omega u^{k+\alpha-1}+\const{a}\int_\Omega |\Delta v|^{\frac{k+\alpha-1}{\alpha-1}}\quad \textrm{for all } t \in (0, \TM).
\end{equation}
We now add to both sides of \eqref{uk_estimate1} the term $\int_\Omega u^k$ and then 
we multiply by $e^t$. Since $e^t \frac{d}{dt}\int_\Omega u^k + e^t \int_\Omega u^k =
\frac{d}{dt}\left(e^t \int_\Omega u^k\right)$ and taking into account bound \eqref{uLap_v_estimate}, an integration over $(0,t)$ provides on $(0,\TM)$
\begin{equation}\label{diff1_u^k}
\begin{split}
&e^t \int_\Omega u^k - \int_\Omega u^k_0 +kb \int_0^t e^s \left(\int_\Omega u^{k+\beta-1}\right)\left(\int_\Omega u^\gamma\right)\,ds
\le -\frac{4(k-1)}{k} \int_0^t e^s  \left(\int_\Omega |\nabla u^{\frac{k}{2}}|^2\right)\,ds\\
&+ \int_0^t e^s \left(\int_\Omega u^k\right)\,ds
+\const{g} \int_0^t e^s \left(\int_\Omega u^{k+\alpha-1}\right)\,ds +\const{a} \int_0^t e^s \left(\int_\Omega|\Delta v|^{\frac{k+\alpha-1}{\alpha-1}}\right)\,ds.
\end{split}
\end{equation}
Let us now address the term $\const{a}  \int_0^t e^s \left( \int_\Omega|\Delta v|^{\frac{k+\alpha-1}{\alpha-1}}\right)ds$.
By applying Lemma \ref{lem:MaxReg} with $q=\frac{k+\alpha-1}{\alpha-1}$, we have for all $t\in(0,\TM)$
\begin{equation*}\label{EstimateLaplV}
\begin{split}
&\const{a} \int_0^t e^s \left(\int_\Omega |\Delta v(\cdot,s)|^{\frac{k+\alpha-1}{\alpha-1}}
\right)ds \le \const{a} C_{P} 2^{\frac{k}{\alpha-1}} \|v_0\|_{L^{\infty}(\Omega)}^{\frac{k+\alpha-1}{\alpha-1}} \int_0^t e^s \left(\int_\Omega (u(\cdot,s))^{\frac{k+\alpha-1}{\alpha-1}}\right)ds + \const{AEI} e^t.
\end{split}
\end{equation*}
From the condition $\alpha \ge 2$ we have that $\frac{k+\alpha-1}{\alpha-1} \le k+\alpha-1$, the Young inequality leads to
\begin{equation*}\label{EstU}
\const{a} C_{P} 2^{\frac{k}{\alpha-1}} \|v_0\|_{L^{\infty}(\Omega)}^{\frac{k+\alpha-1}{\alpha-1}} \int_\Omega u^{\frac{k+\alpha-1}{\alpha-1}} \le 
\const{a} C_{P} 2^{\frac{k}{\alpha-1}} \|v_0\|_{L^{\infty}(\Omega)}^{\frac{k+\alpha-1}{\alpha-1}}
\int_\Omega u^{k+\alpha-1}+\const{p} \quad  \mbox{on } (0, \TM).
\end{equation*}
(For the limiting case $\alpha=2$, the constant $\const{p}$ may naturally be set to $0$.)
Combining the previous estimate with \eqref{diff1_u^k} then yields
\begin{equation}\label{diff2_u^k}
\begin{split}
&e^t \int_\Omega u^k - \int_\Omega u^k_0 +kb \int_0^t e^s \left(\int_\Omega u^{k+\beta-1}\right)\left(\int_\Omega u^\gamma\right)\,ds
\le -\frac{4(k-1)}{k} \int_0^t e^s  \left(\int_\Omega |\nabla u^{\frac{k}{2}}|^2\right)\,ds\\
&+ \int_0^t e^s \left(\int_\Omega u^k\right)\,ds
+\const{g11} \int_0^t e^s \left(\int_\Omega u^{k+\alpha-1}\right)\,ds +\const{p} e^t \quad \textrm{for all } t \in (0,\TM).
\end{split}
\end{equation}
Here, we turn our attention to the third integral on the right-hand side of \eqref{diff2_u^k}. Observing that
$\int_\Omega u^{k+\alpha-1}=\|u^{\frac{k}{2}}\|_{L^{\frac{2(k+\alpha-1)}{k}}(\Omega)}^{\frac{2(k+\alpha-1)}{k}}$, we aim to exploit Lemma \ref{GagliardoIneqLemma} with $\varphi:=u^{\frac{k}{2}}$ and appropriately chosen $q$ and $r$. 
From the definition of $k_0$ and condition \eqref{alpha_large}, for any $k>k_0$ it is feasible to set 
\begin{equation}\label{kPrimo}
k':=\frac{k+\beta+\gamma-1}{2},
\end{equation}
which fulfills 
\begin{equation*}\label{k'condi}
\max\left\{\gamma,\ \frac{k}{2},\ \frac{n(\alpha-1)}{2}\right\}<k'<k+\alpha-1.
\end{equation*}
In this way, for 
\begin{equation*}
q:=\frac{2(k+\alpha-1)}{k},\ r:=\frac{2k'}{k}
\end{equation*}
straightforward calculations reveal $1\le r<q<p$ and $\frac{q}{r}<\frac{2}{r}+1-\frac{2}{p}$. Therefore we infer from \eqref{Bian} that for all $t \in (0,\TM)$
\begin{equation}\label{u^k+alpha-1}
\const{g11} \int_\Omega u^{k+\alpha-1}= \const{g11} \|u^{\frac{k}{2}}\|_{L^{\frac{2(k+\alpha-1)}{k}}(\Omega)}^{\frac{2(k+\alpha-1)}{k}}
\le \frac{2(k-1)}{k}\int_\Omega |\nabla u^{\frac{k}{2}}|^2+ \const{g11} \int_\Omega u^k+\const{i}\left(\int_\Omega u^{k'}\right)^{\frac{\delta}{r}}. 
\end{equation}
Here, the interpolation inequality (see \cite[page 93]{Brezis}) leads on $(0, \TM)$ to 
\begin{equation}\label{u^kp}
\begin{split}
\left(\int_\Omega u^{k'}\right)^{\frac{\delta}{r}}=\|u\|_{L^{k'}(\Omega)}^{b_1}
&\le \|u\|_{L^\gamma(\Omega)}^{a_1b_1} \|u\|_{L^{k+\beta-1}(\Omega)}^{(1-a_1)b_1}\\
&=\left(\|u\|_{L^\gamma(\Omega)}^\gamma \|u\|_{L^{k+\beta-1}(\Omega)}^{k+\beta-1}\right)^{\frac{a_1 b_1}{\gamma}}\|u\|_{L^{k+\beta-1}(\Omega)}
^{\left[1-a_1-\frac{a_1(k+\beta-1)}{\gamma}\right]b_1},
\end{split}
\end{equation}
where 
\begin{equation}\label{DefinitionB1A1}
b_1=b_1(q):=\frac{k'\delta(q)}{r}=\frac{k'\delta}{r}, \quad a_1:=\frac{\frac{1}{k'}-\frac{1}{k+\beta-1}}{\frac{1}{\gamma}-\frac{1}{k+\beta-1}}\in (0,1).
\end{equation}
Recalling the expression of $k'$ in \eqref{kPrimo} and the range of $\gamma$ in \eqref{alpha_large}, some calculations provide 
\begin{equation*}
\left[1-a_1-\frac{a_1(k+\beta-1)}{\gamma}\right]b_1=0 \quad \textrm{and} \quad 
\frac{a_1 b_1}{\gamma}<1. 
\end{equation*}
Consequently, applying Young's inequality, relation \eqref{u^kp} reads
\begin{equation*}
\const{i}\left(\int_\Omega u^{k'}\right)^{\frac{\delta}{r}} \le \const{i}\left(\|u\|_{L^\gamma(\Omega)}^\gamma \|u\|_{L^{k+\beta-1}(\Omega)}^{k+\beta-1}\right)
^{\frac{a_1 b_1}{\gamma}}
\le kb\left(\int_\Omega u^{k+\beta-1}\right)\left(\int_\Omega u^\gamma\right)+\const{j} \quad  \mbox{for\ all}\ t \in (0, \TM),
\end{equation*}
which in conjunction with \eqref{u^k+alpha-1} implies for all $t \in (0, \TM)$, 
\begin{equation}\label{u^k+alpha-1_2}
\const{g11} \int_\Omega u^{k+\alpha-1} \le \frac{2(k-1)}{k}\int_\Omega |\nabla u^{\frac{k}{2}}|^2+ \const{g11} \int_\Omega u^k
+kb\left(\int_\Omega u^{k+\beta-1}\right)\left(\int_\Omega u^\gamma\right)+\const{j}.
\end{equation}
Now we focus on the second integral at the right-hand side: the Gagliardo--Nirenberg inequality and \eqref{MassBounded} produce for
\begin{equation*}
\theta_1:=\frac{\frac{k}{2}-\frac{1}{2}}{\frac{k}{2}+\frac{1}{n}-\frac{1}{2}} \in (0,1)
\end{equation*}
and all $\hat{c}>0$, this bound on $(0,\TM)$:
\begin{equation*}
\hat{c} \int_\Omega u^k= \hat{c}
\|u^\frac{k}{2}\|_{L^2(\Omega)}^2
\le \const{l} \|\nabla u^{\frac{k}{2}}\|_{L^2(\Omega)}^{2\theta_1}\|u^\frac{k}{2}\|_{L^\frac{2}{k}(\Omega)}^{2(1-\theta_1)}
+\const{l}\|u^\frac{k}{2}\|_{L^\frac{2}{k}(\Omega)}^2
\le \const{m}\left(\int_\Omega |\nabla u^{\frac{k}{2}}|^2\right)^{\theta_1}+\const{m}.
\end{equation*}
In turn, we have from the Young inequality that 
\begin{equation}\label{Young}
\begin{split}
\hat{c} \int_\Omega u^k &\le \frac{2(k-1)}{k}\int_\Omega |\nabla u^{\frac{k}{2}}|^2+\const{n} \quad  \mbox{for\ all}\ t \in (0, \TM). 
\end{split}
\end{equation}
By inserting estimate \eqref{u^k+alpha-1_2} into \eqref{diff2_u^k} and taking into account bound \eqref{Young}, we arrive at
\begin{equation*}
\begin{split}
e^t \int_\Omega u^k \leq \const{gs} e^t + \const{gs1}  \quad \textrm{on } (0, \TM),
\end{split}
\end{equation*}
so the claim is proved. 
\end{proof}
\end{lemma}
Let us turn now our attention to the case $1\leq \alpha < 2$.
\begin{lemma}\label{CaseAlphaMinore2Lemma}
Assume that $\alpha, \beta, \gamma \ge 1$ satisfy 
\begin{equation}\label{beta_large}
1 \le \alpha <2, \quad \beta \geq \alpha \quad \mbox{and}\quad \beta + \gamma>\frac{n}{2} +2. 
\end{equation}
Then $u\in L^\infty((0,\TM);L^{k}(\Omega))$, for all $k>k_0$.
\begin{proof}
By following the same argument of Lemma \ref{CaseAlphaMaggiore2Lemma}, we derive for all $k>k_0$ and 
\begin{equation}\label{uk_estimate2}
\begin{split}
\frac{d}{dt}\int_\Omega u^k &=-\frac{4(k-1)}{k}\int_\Omega |\nabla u^{\frac{k}{2}}|^2-(k-1)\chi\int_\Omega u^k\Delta v +ka\int_\Omega u^{k+\alpha-1}-kb\left(\int_\Omega u^{k+\beta-1}\right)\left(\int_\Omega u^\gamma\right)\\
& \leq -\frac{4(k-1)}{k}\int_\Omega |\nabla u^{\frac{k}{2}}|^2 + \int_\Omega u^{k+1}+\const{b}\int_\Omega |\Delta v|^{k+1} +ka\int_\Omega u^{k+\alpha-1}\\ 
&\quad -kb\left(\int_\Omega u^{k+\beta-1}\right)\left(\int_\Omega u^\gamma\right) \quad \textrm{for all } t \in (0,\TM),
\end{split}
\end{equation}
where we used Young's inequality thanks to $\alpha\geq 1$. From the condition $\alpha<2$, an application of Young's inequality to the fourth integral at the right-hand side of \eqref{uk_estimate2} gives 
\begin{equation}\label{Young14}
ka\int_\Omega u^{k+\alpha-1} \le \int_\Omega u^{k+1}+\const{w} \quad  \mbox{on } (0, \TM).
\end{equation}
Combining estimates \eqref{Young14} with bound \eqref{uk_estimate2}, we have for all $t \in (0, \TM)$, 
\begin{equation}\label{uk_estimate3}
\frac{d}{dt}\int_\Omega u^k +kb\left(\int_\Omega u^{k+\beta-1}\right)\left(\int_\Omega u^\gamma\right) \leq -\frac{4(k-1)}{k}\int_\Omega |\nabla u^{\frac{k}{2}}|^2 + 2 \int_\Omega u^{k+1}+\const{b}\int_\Omega |\Delta v|^{k+1}+\const{w}.
\end{equation}
As to the term $\int_\Omega |\Delta v|^{k+1}$ in expression \eqref{uk_estimate3},
by exploiting in this circumstance Lemma \ref{lem:MaxReg} with $q=k+1$, we obtain 
for all $t\in(0,\TM)$
\begin{equation}\label{Estimate1LaplV}
\begin{split}
\const{b} \int_0^t e^s \left(\int_\Omega |\Delta v(\cdot,s)|^{k+1}
\right)ds \le \const{b} C_{P} 2^k \|v_0\|_{L^{\infty}(\Omega)}^{k+1} \int_0^t e^s \left(\int_\Omega (u(\cdot,s))^{k+1}\right) + \const{A12} e^t.
\end{split}
\end{equation}
On the other hand, by adding $\int_\Omega u^k$ at both sides of estimate \eqref{uk_estimate3}, by multiplying what obtained by $e^t$ and considering bound \eqref{Estimate1LaplV}, a subsequent integration over $(0,t)$ yields on $(0,\TM)$
\begin{equation}\label{diff_u^k 143}
\begin{split}
&e^t\int_\Omega u^k-\int_\Omega u_0^k +kb \int_0^t e^s \left(\int_\Omega u^{k+\beta-1}\right)\left(\int_\Omega u^\gamma\right)\,ds\\
&\leq -\frac{4(k-1)}{k} \int_0^t e^s \left(\int_\Omega |\nabla u^{\frac{k}{2}}|^2\right)\,ds 
+ \const{AG} \int_0^t e^s \left(\int_\Omega u^{k+1}\right)\,ds + \int_0^t e^s \left(\int_\Omega u^k\right)\,ds +\const{w1} e^t.
\end{split}
\end{equation}
Now let us focus on the second integral on the right-hand side of \eqref{diff_u^k 143}. Since $\int_\Omega u^{k+1}=\|u^{\frac{k}{2}}\|_{L^{\frac{2(k+1)}{k}}(\Omega)}^{\frac{2(k+1)}{k}}$, we proceed to apply Lemma \ref{GagliardoIneqLemma} with $\varphi:=u^{\frac{k}{2}}$ and suitable $q$ and $r$. In the specific, for any $k>k_0$ and $k'$ given in \eqref{kPrimo}, it is possible to check that 
\begin{equation*}\label{k'condi2}
\max\left\{\gamma, \frac{k}{2}, \frac{n}{2}\right\}<k'<k+1.
\end{equation*} 
In this way, letting
\begin{equation*}
q:=\frac{2(k+1)}{k},\ r:=\frac{2k'}{k}
\end{equation*}
we can establish that $1 \le r<q<p$ and $\frac{q}{r}<\frac{2}{r}+1-\frac{2}{p}$. Consequently,  we deduce from \eqref{Bian} that 
\begin{equation}\label{u^k+1}
\const{AG} \|u^{\frac{k}{2}}\|_{L^{\frac{2(k+1)}{k}}(\Omega)}^{\frac{2(k+1)}{k}} \le \frac{2(k-1)}{k}\int_\Omega |\nabla u^{\frac{k}{2}}|^2+ \const{AG} \int_\Omega u^k+\const{y}\left(\int_\Omega u^{k'}\right)^{\frac{\delta}{r}}\quad  \mbox{for\ all}\ t \in (0, \TM).
\end{equation}
Now an application of the interpolation inequality yields for all $t \in (0, \TM)$, 
\begin{equation*}
\begin{split}
\left(\int_\Omega u^{k'}\right)^{\frac{\delta}{r}}=\|u\|_{L^{k'}(\Omega)}^{b_2} 
&\le \|u\|_{L^\gamma(\Omega)}^{a_2b_2} \|u\|_{L^{k+\beta-1}(\Omega)}^{(1-a_2)b_2}\\
&=\left(\|u\|_{L^\gamma(\Omega)}^\gamma \|u\|_{L^{k+\beta-1}(\Omega)}^{k+\beta-1}\right)^{\frac{a_2b_2}{\gamma}}\|u\|_{L^{k+\beta-1}(\Omega)}
^{\left[1-a_2-\frac{a_2(k+\beta-1)}{\gamma}\right]b_2},
\end{split}
\end{equation*}
where 
\begin{equation*}
b_2=b_2(q):=\frac{k'\delta(q)}{r}=\frac{k'\delta}{r}, \quad 
a_2:=\frac{\frac{1}{k'}-\frac{1}{k+\beta-1}}{\frac{1}{\gamma}-\frac{1}{k+\beta-1}} \in (0,1).
\end{equation*}
(A comparison between the couple $(a_2,b_2)$ above and $(a_1,b_1)$ in \eqref{DefinitionB1A1} shows that $a_1=a_2$, whereas $b_i$, $i=1,2$ depends on $q$.)
From straightforward calculations and the condition \eqref{beta_large}, we note 
that 
\begin{equation*}
\left[1-a_2-\frac{a_2(k+\beta-1)}{\gamma}\right]b_2=0 \quad \textrm{and} \quad \frac{a_2b_2}{\gamma}<1. 
\end{equation*}
Subsequently, we can exploit the Young inequality entailing 
\begin{equation*}
\const{y}\left(\int_\Omega u^{k'}\right)^{\frac{\delta}{r}}\le\const{y}\left(\|u\|_{L^\gamma(\Omega)}^\gamma\|u\|_{L^{k+\beta-1}(\Omega)}^{k+\beta-1}\right)^{\frac{a_2b_2}{\gamma}}
\le kb\left(\int_\Omega u^{k+\beta-1}\right)\left(\int_\Omega u^\gamma\right)+\const{z} \quad  \mbox{on }\;(0, \TM).
\end{equation*}
This in conjunction with \eqref{u^k+1} implies that 
$t \in (0, \TM)$, 
\begin{equation}\label{u^k+1_2}
\const{AG} \int_\Omega u^{k+1}\le \frac{2(k-1)}{k} \int_\Omega |\nabla u^{\frac{k}{2}}|^2+ \const{AG} \int_\Omega u^k +
k b \left(\int_\Omega u^{k+\beta-1}\right)\left(\int_\Omega u^\gamma\right)
+\const{z}. 
\end{equation}
By rearranging bound \eqref{diff_u^k 143} by virtue of estimates \eqref{u^k+1_2} and \eqref{Young}, it is provided
\begin{equation*}
\int_\Omega u^k \leq \const{ad} + \const{ad1} e^{-t} \quad \textrm{on }\, (0,\TM),
\end{equation*}
which gives the statement.
 \end{proof}
\end{lemma}
\subsection{The case $f(u)=-a u^\alpha+bu^\beta \int_\Omega u^\gamma$}
\begin{lemma}\label{LemmaProbInv}
Assume $a, b >0$ and $\alpha, \beta, \gamma \geq 1$ satisfy either
\begin{equation}\label{CondY}
\alpha >2, \quad \beta + \gamma < \alpha,
\end{equation} 
or 
\begin{equation}\label{Cond1Y}
\alpha >2, \quad \beta + \gamma = \alpha \quad \textrm{and} \quad a > b |\Omega|.
\end{equation} 
Then $u \in L^{\infty}((0,\TM); L^k(\Omega))$ for all $k>k_0$.
\begin{proof}
By reasoning as in the proof of Lemma \ref{CaseAlphaMaggiore2Lemma}, 
for all $k>k_0$ we can write 
\begin{equation}\label{E_inv1}
\begin{split}
\frac{d}{dt}\int_\Omega u^k &=-\frac{4(k-1)}{k}\int_\Omega |\nabla u^{\frac{k}{2}}|^2-(k-1)\chi\int_\Omega u^k\Delta v + kb \int_\Omega u^{k+\beta-1} \int_\Omega u^\gamma -ka\int_\Omega u^{k+\alpha-1}\\
& \leq (\epsilon-ka) \int_\Omega u^{k+\alpha-1} 
+ \const{GS} \int_\Omega |\Delta v|^{\frac{k+\alpha-1}{\alpha-1}}
+ kb \int_\Omega u^{k+\beta-1} \int_\Omega u^\gamma \quad \textrm{for all } t \in (0,\TM),
\end{split}
\end{equation}
where again we have applied the Young inequality (recall $\alpha >1$). 
Adding $\int_\Omega u^k$ to both sides of estimate \eqref{E_inv1} and multiplying the resulting expression by $e^t$, integration over $(0,t)$ then yields, for all
$t \in (0, \TM)$,
\begin{equation}\label{E_inv5}
\begin{split}
&e^t\int_\Omega u^k-\int_\Omega u_0^k \leq 
(\epsilon - k a)  \int_0^t e^s \left(\int_\Omega u^{k+\alpha-1}\right)\, ds
+ \int_0^t e^s \left(\int_\Omega u^k\right)\,ds\\
&+ \const{GS} \int_0^t e^s \left(\int_\Omega |\Delta v|^{\frac{k+\alpha-1}{\alpha-1}}\right)\, ds + kb \int_0^t e^s \left(\int_\Omega u^{k+\beta-1} \int_\Omega u^\gamma\right)\, ds.
\end{split}
\end{equation}
By applying Lemma \ref{lem:MaxReg} with $q=\frac{k+\alpha-1}{\alpha-1}$ and Young's inequality since $\alpha >2$, we obtain on $(0,\TM)$
\begin{equation}\label{E_inv6}
\begin{split}
\const{GS} \int_0^t e^s \left(\int_\Omega |\Delta v|^{\frac{k+\alpha-1}{\alpha-1}}\right)\, ds &\leq \const{AIi} e^t + \const{GS} C_P 2^{\frac{k}{\alpha-1}} \|v_0\|_{L^{\infty}(\Omega)}^{\frac{k+\alpha-1}{\alpha-1}}  \int_0^t e^s \left(\int_\Omega u^{\frac{k+\alpha-1}{\alpha-1}}\right)\, ds\\
& \leq \const{AIiB} e^t + \epsilon \int_0^t e^s \left(\int_\Omega u^{k+\alpha-1} \right)\, ds.
\end{split}
\end{equation}
A further application of Young's inequality yields 
\begin{equation}\label{E_inv7}
\int_0^t e^s \left(\int_\Omega u^k\right)\,ds \leq \epsilon \int_0^t e^s \left(\int_\Omega u^{k+\alpha-1}\right)\,ds + \const{AEi} e^t \quad \textrm{for all } t \in (0,\TM).
\end{equation}
By virtue of $\beta<\beta+\gamma<\alpha$ and $k>k_0$, a combination of H\"{o}lder's and Young's inequalities leads for all $t \in (0,\TM)$ to  
\begin{equation}\label{E_inv3}
kb \int_\Omega u^{k+\beta-1} \int_\Omega u^\gamma
\leq 
\begin{cases}
k b |\Omega|^{\frac{k+2\alpha-\beta-\gamma-1}{k+\alpha-1}} \left( \displaystyle \int_\Omega u^{k+\alpha-1}\right)^{\frac{k+\beta+\gamma-1}{k+\alpha-1}} \leq \epsilon  \displaystyle \int_\Omega u^{k+\alpha-1} + \const{AI} & \textrm{if $\beta + \gamma < \alpha$},\\
k b |\Omega|^{\frac{k+2\alpha-\beta-\gamma-1}{k+\alpha-1}} \left(\displaystyle \int_\Omega u^{k+\alpha-1}\right)^{\frac{k+\beta+\gamma-1}{k+\alpha-1}} = k b |\Omega| \displaystyle \int_\Omega u^{k+\alpha-1} & \textrm{if $\beta + \gamma = \alpha$}.
\end{cases}
\end{equation}
By plugging estimates \eqref{E_inv6}, \eqref{E_inv7} and \eqref{E_inv3} into \eqref{E_inv5}, we derive on $(0,\TM)$
\begin{equation*}\label{E_inv8}
e^t\int_\Omega u^k\leq 
\begin{cases}
\const{HE} + \const{AIcd} e^t + (\epsilon-ka) \displaystyle \int_0^t e^s \left(\int_\Omega u^{k+\alpha-1}\right)\, ds & \textrm{if $\beta + \gamma < \alpha$}\\
\const{HE} + \const{AIcd} e^t + (\epsilon + kb |\Omega|-ka) \displaystyle \int_0^t e^s \left(\int_\Omega u^{k+\alpha-1}\right)\, ds & \textrm{if $\beta + \gamma = \alpha$},
\end{cases}
\end{equation*}
which implies, by virtue of the arbitrariness of $\epsilon$ and assumptions \eqref{CondY} and \eqref{Cond1Y}, 
\begin{equation*}
e^t \int_\Omega u^k \leq \const{HE} + \const{AIcd} e^t \quad \textrm{on }\, (0,\TM),
\end{equation*}
concluding that $u \in L^{\infty}((0,\TM); L^k(\Omega))$ for all $k>k_0$.
\end{proof}
\end{lemma}
In the following lemma, we set $\alpha = 2$ and $\beta + \gamma = \alpha$. Under the conditions $\beta \ge 1$ and $\gamma \ge 1$, this uniquely determines $\beta = \gamma = 1$.
\begin{lemma}\label{LemmaProbInv1}
Assume $\alpha=2$, and  $\beta = \gamma =1$, $k>k_0$. Moreover, let $C_P=C_P(n, \Omega, k+1)$ being $C_P$ the constant in Lemma \ref{lem:MaxReg}.
If 
\begin{equation}\label{Cond0Y}
a > b |\Omega| + \frac{(k-1)}{k} (C_P 2^{k})^{\frac{1}{k+1}} \|\chi v_0\|_{L^{\infty}(\Omega)},
\end{equation} 
then the same statement of Lemma \ref{LemmaProbInv} holds.
\begin{proof}
By taking $\alpha=2$ and $\beta=\gamma=1$ in bounds \eqref{E_inv1} and \eqref{E_inv3} , we obtain 
for all $t \in (0, \TM)$ 
\begin{equation*}\label{Ec_inv1}
\begin{split}
\frac{d}{dt}\int_\Omega u^k &\leq \left(\frac{k(k-1)}{(k+1)} (C_P 2^{k})^{\frac{1}{k+1}} \|\chi v_0\|_{L^{\infty}(\Omega)} - ka \right) \int_\Omega u^{k+1} \\
&\quad + \frac{(k-1)\chi}{k+1}\left(C_P 2^k\right)^{-\frac{k}{k+1}} \|v_0\|_{L^{\infty}(\Omega)}^{-k}  
\int_\Omega |\Delta v|^{k+1}+ kb \int_\Omega u^k \int_\Omega u\\
& \leq  \left(\frac{k(k-1)}{(k+1)} (C_P 2^{k})^{\frac{1}{k+1}} \|\chi v_0\|_{L^{\infty}(\Omega)} + kb |\Omega| - ka \right) \int_\Omega u^{k+1} \\
&\quad + \frac{(k-1)\chi}{k+1}\left(C_P 2^k\right)^{-\frac{k}{k+1}} \|v_0\|_{L^{\infty}(\Omega)}^{-k}  
\int_\Omega |\Delta v|^{k+1}.
\end{split}
\end{equation*}
Following the argument used in the proof of Lemma~\ref{LemmaProbInv}, we obtain for all
$t \in (0, \TM)$
\begin{equation}\label{Ec_inv5}
\begin{split}
&e^t\int_\Omega u^k-\int_\Omega u_0^k \leq 
\left(\frac{k(k-1)}{(k+1)} (C_P 2^{k})^{\frac{1}{k+1}} \|\chi v_0\|_{L^{\infty}(\Omega)} + kb |\Omega| - ka \right) \int_0^t e^s \left(\int_\Omega u^{k+1}\right)\, ds
\\
&+  \frac{(k-1)\chi}{k+1}\left(C_P 2^k\right)^{-\frac{k}{k+1}} \|v_0\|_{L^{\infty}(\Omega)}^{-k} \int_0^t e^s \left(\int_\Omega |\Delta v|^{k+1}\right)\, ds +  \int_0^t e^s \left(\int_\Omega u^k\right)\,ds.
\end{split}
\end{equation}
Now, an application of Lemma \ref{lem:MaxReg} with $q=k+1$ entails on $(0,\TM)$
\begin{equation}\label{Ec_inv6}
\begin{split}
&\frac{(k-1)\chi}{k+1}\left(C_P 2^k\right)^{-\frac{k}{k+1}} \|v_0\|_{L^{\infty}(\Omega)}^{-k}    \int_0^t e^s \left(\int_\Omega |\Delta v|^{k+1}\right)\, ds\\ &\leq \const{AIi} e^t + \frac{(k-1)}{k+1}\left(C_P 2^k\right)^{\frac{1}{k+1}} \|\chi v_0\|_{L^{\infty}(\Omega)}  
\int_0^t e^s \left(\int_\Omega u^{k+1}\right)\, ds.
\end{split}
\end{equation}
By plugging estimates \eqref{Ec_inv6} and \eqref{E_inv7} into \eqref{Ec_inv5}, we have on $(0,\TM)$
\begin{equation*}\label{Ec_inv8}
\begin{split}
&e^t\int_\Omega u^k\leq \const{HE} + \const{AIcd} e^t + k
\left[\frac{(k-1)}{(k+1)} (C_P 2^{k})^{\frac{1}{k+1}} \|\chi v_0\|_{L^{\infty}(\Omega)}
+ b |\Omega| + \epsilon \right.\\
&\left. + \frac{(k-1)}{k(k+1)}\left(C_P 2^k\right)^{\frac{1}{k+1}} \|\chi v_0\|_{L^{\infty}(\Omega)} - a \right]  \int_0^t e^s \left(\int_\Omega u^{k+1}\right)\, ds.
\end{split}
\end{equation*}
Finally, the assumption \eqref{Cond0Y} made on $a$ allows us to achieve the claim.
\end{proof}
\end{lemma}
\section{Proof of Theorem \ref{theoremlocal}}
The claims concerning \ref{a} and \ref{b} for the case where $f(u)=au^\alpha-bu^\beta \int_\Omega u^\gamma$, are consequence of Lemma \ref{CaseAlphaMaggiore2Lemma} and Lemma \ref{CaseAlphaMinore2Lemma} combined with Lemma \ref{ExtensionLemma}. 

For $f(u)=-au^\alpha + bu^\beta \int_\Omega u^\gamma$ the cases \ref{(a)} and \ref{(b)} follow by reasoning as before, invoking Lemma \ref{LemmaProbInv} and Lemma \ref{ExtensionLemma}. Finally, let $\mathcal{C}_P= \mathcal{C}_P(n, \Omega)=
\frac{n+1}{n+2} (C_P (n, \Omega, n+3) \, 2^{n+2})^{\frac{1}{n+3}}$, being  
$C_P$ introduced in Lemma \ref{lem:MaxReg}. Under the assumption on $a$ in case \ref{(c)}, by continuity arguments also involving the constant $C_P$ itself, there is $k>n+2$ for which \eqref{Cond0Y} holds, so that by exploiting Lemma \ref{LemmaProbInv1} we conclude.
\newpage
\subsection*{\bf\textit{Conflict of interest disclosure}}
No potential conflict of interest was reported by the authors
\subsubsection*{\bf\textit{\quad Acknowledgments}}
The authors are members of the Gruppo Nazionale per l'Analisi Matematica, la Probabilit\`a e le loro Applicazioni (GNAMPA) of the Istituto Nazionale di Alta Matematica (INdAM), and are partially supported by the research project {\em Partial Differential Equations and their role in understanding natural phenomena} (2023, CUP F23C25000080007), funded by  \href{https://www.fondazionedisardegna.it/}{Fondazione di Sardegna}. 
RDF acknowledges financial support by PNRR e.INS Ecosystem of Innovation
for Next Generation Sardinia (CUP F53C22000430001, codice MUR ECS0000038).
SF and GV are also supported by MIUR (Italian Ministry of Education, University and Research) Prin 2022 \emph{Nonlinear differential problems with applications to real phenomena} (Grant Number: 2022ZXZTN2). 

\end{document}